\documentclass[12pt, twoside, eqno]{article}
\usepackage{latexsym}
\usepackage{amssymb}
\usepackage{amsfonts}
\begin{document}

\noindent {\bf \large Decomposition of intra-regular
$po$-$\Gamma$-semigroups into simple components}\bigskip

\medskip

\noindent{\bf Niovi Kehayopulu}\bigskip

\noindent November 29, 2015\bigskip{\small

\noindent{\bf Abstract} We keep the definition of intra-regularity 
(left regularity) of $po$-$\Gamma$-semigroups introduced in [6] which 
is absolutely necessary for the investigation. Being able to describe 
the form of the elements of the principal filter by using this 
definition, we study the decomposition of an intra-regular 
$po$-$\Gamma$-semigroup into simple components.
Then we prove that a $po$-$\Gamma$-semigroup $M$ is intra-regular and 
the ideals of $M$ form a chain if and only if $M$ is a chain of 
simple semigroups. Moreover, a $po$-$\Gamma$-semigroup $M$ is 
intra-regular and the ideals of $M$ form a chain if and only if the 
ideals of $M$ are prime. Finally, for an intra-regular 
$po$-$\Gamma$-semigroup $M$, the set $\{(x)_{\cal N} \mid x\in M\}$ 
coincides with the set of all maximal simple subsemigroups of $M$. A 
decomposition of left regular and left duo $po$-$\Gamma$-semigroup 
into left simple components has been also given. \medskip

\medskip

 }
\section{Introduction and prerequisites} An ordered 
$\Gamma$-semigroup (shortly $po$-$\Gamma$-semigroup) is a nonempty 
set $M$ with a set of binary operations $\Gamma$ on $M$ and an order 
relation on $M$ such that $a\le b$ implies $a\gamma c\le b\gamma c$ 
and $c\gamma a\le c\gamma b$ for every $c\in M$ and every 
$\gamma\in\Gamma$. For a $po$-$\Gamma$-semigroup $M$ and a subset $H$ 
of $M$ we denote by $(H]$ the subset of $M$ defined by $(H]=\{t\in M 
\mid t\le a \mbox { for some } a\in H\}$. We have $M=(M]$, and for 
any two subsets $A$, $B$ of $M$, we have $A\subseteq (A]$; if $A$ is 
a right (or left) ideal of $M$, then $A=(A]$; $(A]\Gamma (B]\subseteq 
(A\Gamma B]$; $A\subseteq B$ implies $(A]\subseteq (B]$; $((A]]=(A]$. 
Let $M$ be a $po$-$\Gamma$-semigroup. A subset $A$ of $M$ is called a 
{\it left} (resp. {\it right}) {\it ideal} of $M$ if (1) $M\Gamma 
A\subseteq A$ and (2) if $a\in A$ and $M\ni b\le a$, then $b\in A$. 
It is called an {\it ideal} of $M$ if it is both a left and right 
ideal of $M$. For an element $a$ of $M$, we denote by $L(a)$, $R(a)$ 
and $I(a)$ the left ideal, right ideal and the ideal of $M$, 
respectively, generated by $a$, and we have $L(a)=(a\cup M\Gamma a]$, 
$R(a)=(a\cup a\Gamma M]$, and $I(a)=(a\cup M\Gamma a\cup a\Gamma 
M\cup M\Gamma a\Gamma M]$. We denote by ${\cal L}$ the equivalence 
relation on $M$ defined by $a{\cal L}b$ if and only if $L(a)=L(b)$, 
by ${\cal R}$ the equivalence relation on $M$ defined by $a{\cal R}b$ 
if and only if $R(a)=R(b)$ and by ${\cal I}$ the equivalence relation 
on $M$ defined by $a{\cal I} b$ if and only if $I(a)=I(b)$. A 
nonempty subset $A$ of $M$ is called a {\it subsemigroup} of $M$ is 
for every $a,b\in A$ and every $\gamma\in\Gamma$ we have $a\gamma 
b\in A$, that is, if $A\Gamma A\subseteq A$. A subsemigroup $F$ of 
$M$ is called a {\it filter} of $M$ if (1) $a,b\in F$ and 
$\gamma\in\Gamma$ such that $a\gamma b\in F$ implies $a\in F$ and 
$b\in F$ and (2) if $a\in F$ and $M\ni c\ge a$, then $c\in F$. An 
equivalence relation $\sigma$ on $M$ is called {\it congruence} if 
$(a,b)\in\sigma$ implies $(a\gamma c,b\gamma c)\in\sigma$ and 
$(c\gamma a,c\gamma b)\in\sigma$ for any $c\in M$ and any $\gamma\in 
\Gamma$. A congruence $\sigma$ on $M$ is called {\it semilattice 
congruence} if $(a\gamma b, b\gamma a)\in \sigma$ and $(a\gamma a, 
a)\in\sigma$ for every $a,b\in M$ and every $\gamma\in\Gamma$. If 
$\sigma$ is a semilattice congruence on $M$, then the $\sigma$-class 
$(a)_\sigma$ of $M$ containing the element $a$ is a subsemigroup of 
$M$ for every $a\in M$. A semilattice congruence $\sigma$ on $M$ is 
called {\it complete} if $a\le b$ implies $(a,a\gamma b)\in\sigma$ 
for every $\gamma\in\Gamma$. We denote by $\cal N$ the relation on 
$M$ defined by $a{\cal N} b$ if and only if the filters of $M$ 
generated by the elements $a$ and $b$ of $M$ coincide. As in 
semigroups, the relation ${\cal N}$ is a semilattice congruence on 
$M$.
So, if $z\in M$ and $\gamma\in\Gamma$, then we have $(z\gamma z,z)\in 
{\cal N}$, $(z\gamma z\gamma z,z\gamma z)\in {\cal N}$, $(z\gamma 
z\gamma z\gamma z,z\gamma z\gamma z)\in {\cal N}$ and so on.
In particular, exactly as in ordered semigroups, the relation $\cal 
N$ is a complete semilattice congruence on $M$. A $\Gamma$-semigroup 
$M$ is called {\it left} (resp. {\it right}) {\it simple} if for 
every left (resp. right) ideal $T$ of $M$ we have $T=M$, that is, if 
$M$ is the only left (resp. right) ideal of $M$. A $\Gamma$-semigroup 
$M$ is called {\it simple} if $M$ is the only ideal of $M$.
A $po$-$\Gamma$-semigroup $M$ is said to be a {\it semilattice of 
simple} (resp. {\it left simple}) {\it semigroups} if there exists a 
semilattice congruence $\sigma$ on $M$ such that the $\sigma$-class 
$(x)_{\sigma}$ is a simple (resp. left simple) subsemigroup of $M$ 
for every $x\in M$. A $po$-$\Gamma$-semigroup $M$ is called a {\it 
chain of simple} (resp. {\it left simple}) {\it semigroups} if there 
exists a semilattice congruence $\sigma$ on $M$ such that 
$(x)_\sigma$ is a simple (resp. left simple) subsemigroup of $M$ for 
every $x\in M$, and for any $x,y\in M$ and any $\gamma\in\Gamma$ we 
have $(x)_{\sigma}=(x\gamma y)_{\sigma}$ or $(y)_{\sigma}=(x\gamma 
y)_{\sigma}$.

Many results on $\Gamma$-semigroups can be obtained from semigroups 
just putting a ``Gamma" in the appropriate place. But there are also 
results for which the transfer is not so easy. In a 
$\Gamma$-semigroup $M$ the filter of $M$ generated by an element $a$ 
of $M$ plays an essential role in the structure, in particular, in 
the decomposition of some types of $\Gamma$-semigroups. So it is 
important to get the form of its elements. The existed definition of 
intra-regularity in the bibliography used to be the following: A 
$po$-$\Gamma$ semigroup $M$ is intra-regular if, by definition, 
$(a\in M\Gamma a\Gamma a\Gamma M]$ for every $a\in M$. With this 
definition is not possible to describe the form of the elements of 
the $N(x)$ $(x\in S)$. To overcome this difficulty, the following new 
concept of intra-regularity has been introduced in [6]: We say that a 
$po$-$\Gamma$-semigroup $M$ is intra-regular if $a\in (M\Gamma 
a\gamma a\Gamma M]$ for  every $a\in M$. Using this definition, we 
first give a structure theorem referring to the decomposition of a 
$po$-$\Gamma$-semigroup into simple components. Then, for a 
$po$-$\Gamma$-semigroup $M$, we prove the following: The ideals of 
$M$ are prime if and only if they form a chain and $M$ is 
intra-regular. $M$ is intra-regular and the ideals of $M$ form a 
chain if and only if $M$ is a chain of simple semigroups. For an 
intra-regular $po$-$\Gamma$-semigroup $M$ the set $\{(x)_{\cal N} 
\mid x\in M\}$ coincides with the set of all maximal simple 
subsemigroups of $M$. In our investigation, we use the fact that if 
the ideals of a $po$-$\Gamma$-semigroup are weakly prime, then they 
form a chain. A $po$-$\Gamma$-semigroup $M$ is called {\it left} 
(resp. {\it right}) {\it duo} if the left (resp. right) ideals of $M$ 
are two-sided. Keeping the new definition of left (right) regularity 
of $po$-$\Gamma$-semigroups introduced in [6], we also give a 
structure theorem related with the decomposition of a 
$po$-$\Gamma$-semigroup $M$ which is left regular and left duo into 
left simple components. If we want to get a result on a 
$po$-$\Gamma$-semigroup, we never work directly on the 
$po$-$\Gamma$-semigroup. Exactly as in the hypersemigroups, we have 
to solve it first for an ordered semigroup and then to be careful to 
define the analogous concepts in case of $po$-$\Gamma$-semigroups (if 
they do not defined directly) and put the ``$\Gamma$" in the 
appropriate place. We never solve the problem directly in 
$po$-$\Gamma$-semigroups.
The results of this paper are based on the corresponding results on 
ordered semigroups considered in [2] and [3], and the aim of writing 
this paper is to show the importance of these new concepts of 
intra-regularity and left (right) regularity considered in [6], in 
the investigation. \section{Main results}
\noindent Let $M$ be a $po$-$\Gamma$-semigroup. A subset $A$ of $M$ 
is called {\it idempotent}, if $A=(A\Gamma A]$. A subset $T$ of $M$ 
is called {\it prime} if $a, b\in M$ and $\gamma\in\Gamma$ such that 
$a\gamma b\in T$ implies $a\in T$ or $b\in T$. The set $T$ is called 
{\it semiprime} if $a\in M$ and $\gamma\in\Gamma$ such that $a\gamma 
a\in T$ implies $a\in T$ [6]. A subset $T$ of $M$ is called {\it 
weakly prime} if the following assertion is satisfied:

If $A$, $B$ are ideals of $M$ such that $A\Gamma B\subseteq T$, then 
$A\subseteq T$ or $B\subseteq T$.\medskip

\noindent For a subset $T$ of $M$, we consider the statements:

(1) $a, b\in M$, $\gamma\in\Gamma$, $a\gamma b\in T$ 
$\Longrightarrow$ $a\in T$ or $b\in T$.

(2) $A, B\subseteq M$, $A\Gamma B\subseteq T$ $\Longrightarrow$ 
$A\subseteq T$ or $B\subseteq T$.\smallskip

\noindent Then $(1)\Rightarrow (2)$. In fact: Let $A, B\subseteq M$, 
$A\Gamma B\subseteq T$, $A\nsubseteq T$ and $b\in B$. Take an element 
$a\in A$ such that $a\not\in T$ and an element $\gamma\in\Gamma$ 
$(\Gamma\not=\emptyset)$. Since $a\gamma b\in A\Gamma B\subseteq T$, 
by (1), we have $a\in T$ or $b\in T$. Since $a\not\in T$, we get 
$b\in T$.\medskip

\noindent We have the following:

(a) If $T$ is a prime subset of $M$, then $T$ is a semiprime subset 
of $M$.

(b) If $T$ is a prime subset of $M$, then $T$ is a weakly prime 
subset of $M$. \medskip

\noindent{\bf Definition 1.} [6] A $po$-$\Gamma$-semigroup $M$ is 
called {\it intra-regular} if$$x\in (M\Gamma x\gamma x\Gamma M]$$for 
every $x\in M$ and every  $\gamma\in\Gamma$.\medskip

\noindent{\bf Proposition 2.} {\it If M is an intra-regular 
$po$-$\Gamma$-semigroup then, for every $x,y\in M$ and every 
$\gamma\in\Gamma$, we have $(M\Gamma x\gamma y\Gamma M]=(M\Gamma 
y\gamma x\Gamma M].$}

\noindent {\bf Proof.} Let $x,y\in M$ and $\gamma\in\Gamma$. Since 
$x\gamma y\in M\Gamma M\subseteq M$ and $M$ is intra-regular, we have 
$x\gamma y\in {\Big(}M\Gamma (x\gamma y)\gamma (x\gamma y)\Gamma 
M{\Big]}\subseteq (M\Gamma y\gamma x\Gamma M].$ Then we have 
\begin{eqnarray*}M\Gamma (x\gamma y)\Gamma M&\subseteq&(M]\Gamma 
(M\Gamma y\gamma x\Gamma M]\Gamma (M]\subseteq (M\Gamma M\Gamma 
y\gamma x\Gamma M\Gamma M]\\&\subseteq& (M\Gamma y\gamma x\Gamma 
M],\end{eqnarray*}from which $(M\Gamma (x\gamma y)\Gamma M]\subseteq 
{\Big (}(M\Gamma y\gamma x\Gamma M]{\Big]}=(M\Gamma y\gamma x\Gamma 
M].$ Since $M$ is intra-regular and $y\gamma x\in M$, by symmetry, we 
get $(M\Gamma y\gamma x\Gamma M]\subseteq (M\Gamma x\gamma y\Gamma 
M]$, thus we have $(M\Gamma x\gamma y\Gamma M]=(M\Gamma y\gamma 
x\Gamma M]$.$\hfill\Box$\medskip

\noindent{\bf Lemma 3.} [6] {\it A $po$-$\Gamma$-semigroup $M$ is 
intra-regular if and only if, for every $x\in M$, we have $$N(x)=
\{y\in M \mid x\in (M\Gamma y\Gamma M]\}.$$}In a similar way as in 
[1] and [4], we can prove the following lemma.\medskip

\noindent{\bf Lemma 4.} {\it If M is a $po$-$\Gamma$-semigroup, then 
${\cal I}\subseteq {\cal N}$ and ${\cal I}\subseteq {\cal 
L}$.}\medskip

\noindent{\bf Lemma 5.} [6] {\it A $po$-$\Gamma$-semigroup $M$ is 
intra-regular if and only if the ideals of M are semiprime}.\\The 
proof of the following lemma is easy.\medskip

\noindent{\bf Lemma 6.} {\it If $M$ is a $po$-$\Gamma$-semigroup, 
then the set $(M\Gamma a\Gamma M]$ (resp. $(M\Gamma a]$) is an ideal 
(resp. left ideal) of $M$, and the set $(a\Gamma M]$ is a right ideal 
of M} for every $a\in M$.\medskip

\noindent{\bf Definition 7.} A $po$-$\Gamma$-semigroup $S$ is said to 
be a {\it semilattice of simple} (resp. {\it left simple}) {\it 
semigroups} if there exists a semilattice congruence $\sigma$ on $M$ 
such that the $\sigma$-class $(x)_{\sigma}$ of $M$ containing $x$ is 
a simple (resp. left simple) subsemigroup of $M$ for every $x\in 
M$.\medskip

\noindent{\bf Theorem 8.} {\it Let M be an $po$-$\Gamma$-semigroup. 
The following are equivalent:

$(1)$ M is intra-regular.

$(2)$ $N(x)=\{y\in M \mid x\in (M\Gamma y\Gamma M]\}$ for every
$x\in M$.

$(3)$ $\cal N=\cal I$.

$(4)$ For every ideal I of M, we have $I=
\bigcup\limits_{x \in I} {(x)_{\cal N}}$.

$(5)$ $(x)_{\cal N}$ is a simple subsemigroup of M for
every $x\in M$.

$(6)$ $M$ is a semilattice of simple semigroups.

$(7)$ Every ideal of M is semiprime.}\medskip

\noindent{\bf Proof.} The implication $(1)\Rightarrow (2)$ follows 
from Lemma 3, the proof of $(3)\Rightarrow (4)$ is similar with the 
corresponding result for semigroups without order in [6], 
$(5)\Rightarrow (6)$ since $\cal N$ is a semilattice congruence on 
$M$ and $(7)\Rightarrow (1)$ by Lemma 5.\smallskip

\noindent$(2)\Longrightarrow (3)$. Let $(a,b)\in {\cal N}$. Since 
$a\in N(a)=N(b)$, by (2), we have $b\in (M\Gamma a\Gamma M]\subseteq 
(a\cup M\Gamma a\cup a\Gamma M\cup M\Gamma a\Gamma M]=I(a)$, so 
$I(b)\subseteq I(a)$. Since $b\in N(a)$, by symmetry, we get 
$I(a)\subseteq I(b)$, so $I(a)=I(b)$, and $(a,b)\in{\cal I}$. Then 
${\cal N}\subseteq {\cal I}$, on the other hand by Lemma 4, we have 
${\cal I}\subseteq {\cal N}$, thus ${\cal I}={\cal N}$. \smallskip

\noindent$(4)\Longrightarrow (5)$. Let $x\in M$. Since $\cal N$ is a 
semilattice congruence on $M$, $(x)_{\cal N}$ is a subsemigroup of 
$M$. Let $I$ be an ideal of $(x)_{\cal N}$. Then $I=(x)_{\cal N}$. 
Indeed: Let $y\in (x)_{\cal N}$. Take an element $z\in I$ and an 
element $\gamma\in\Gamma$ $(\Gamma\not=\emptyset)$. Since
($M\Gamma z\gamma z\gamma z\Gamma M]$ is an ideal of $M$ (see Lemma 
6), by hypothesis, we have $(M\Gamma z\gamma z\gamma z\Gamma M]=
\bigcup\limits_{t \in (M\Gamma z\gamma z\gamma z\Gamma M]}
{(t)_{\cal N}}$. Since

$y\in (x)_{\cal N}=(z)_{\cal N}=(z\gamma z\gamma z\gamma z\gamma 
z)_{\cal N}\subseteq (M\Gamma z\gamma z\gamma z\Gamma M],$ we have 
$$y\le a\delta z\gamma z\gamma
z\xi b=(a\delta z)\gamma z\gamma (z\xi b)\; \mbox { for some }
a,b\in M, \;\delta,\xi\in\Gamma.$$Using the fact that ${\cal N}$ is a 
complete semilattice congruence on $M$, in a similar way as in [5], 
we prove that $a\delta z\in (x)_{\cal N}$ and $z\xi b\in (x)_{\cal 
N}$. Then, since $I$ is an ideal of $(x)_{\cal N}$ and $z\in I$, we 
have
we have $(a\delta z)\gamma z\gamma (z\xi b)\in (x)_{\cal N}\Gamma
I\Gamma (x)_{\cal N}\subseteq I$, and $y\in I$. Hence $(x)_{\cal 
N}\subseteq I$, and so $I=(x)_{\cal N}$.\smallskip

\noindent$(6)\Longrightarrow (7)$. Let $\sigma$ be a semilattice 
congruence on $M$ such that $(x)_\sigma$ is a simple subsemigroup of 
$M$ for every $x\in M$. Let $I$ be an ideal of $M$,
$x\in M$ and $\gamma\in\Gamma$ such that $x\gamma x\in I$. The set 
$I\cap (x)_{\sigma}$ is an ideal of $(x)_\sigma$. Indeed: Taking into 
account the proof of the implication $(6)\Rightarrow (7)$ in [5], it 
is enough to prove the following: Let $a\in I\cap (x)_{\sigma}$ and 
$(x)_{\sigma}\ni b\le a$, then $b\in I\cap (x)_{\sigma}$. Since $M\in 
b\le a\in I$ and $I$ is an ideal of $M$, we have $b\in I$, then $b\in 
I\cap (x)_{\sigma}$. Since $(x)_\sigma$ is a simple subsemigroup of 
$M$, we have $I\cap (x)_\sigma=(x)_\sigma$, then $x\in I$. Thus $M$ 
is semiprime. $\hfill\Box$

A subset $A$ of a $\Gamma$-semigroup $M$ is called {\it idempotent} 
if $A=(A\Gamma A]$.\medskip

\noindent {\bf Lemma 9.} {\it Let M be a $po$-$\Gamma$-semigroup. The 
ideals of M are idempotent if and only if for any ideals $A, B$ of M, 
we have $A\cap B=(A\Gamma B]$.}\medskip

\noindent{\bf Proof.} $\Longrightarrow$. Let $A, B$ be ideals of $M$. 
Then $(A\Gamma B]\subseteq (A\Gamma M]\subseteq (A]=A$ and $(A\Gamma 
B]\subseteq (M\Gamma B]\subseteq (B]=B$, thus $(A\Gamma B]\subseteq 
A\cap B$.
On the other hand, $A\cap B$ is an ideal of $M$. Indeed: Take an 
element $a\in A$, an element $b\in B$ and an element 
$\gamma\in\Gamma$ $(A, B, \Gamma\not=\emptyset)$. Then $a\gamma b\in 
A\Gamma B\subseteq A\Gamma M\subseteq A$ and $a\gamma b\in A\Gamma 
B\subseteq M\Gamma B\subseteq B$, so $a\gamma b\in A\cap B$, and 
$\emptyset\not=A\cap B\subseteq M$. We also have $(A\cap B)\Gamma 
M\subseteq A\Gamma M\subseteq A$, $M\Gamma (A\cap B)\subseteq M\Gamma 
B\subseteq B$, and if $x\in A\cap B$ and $M\ni y\le x$ then, since 
$x\in A$ we have $y\in A$ and since $x\in B$ we have $y\in B$, so 
$y\in A\cap B$. Since $A\cap B$ is an ideal of $M$, by hypothesis, we 
have $A\cap B={\Big(}(A\cap B)\Gamma (A\cap B){\Big]}\subseteq 
(A\Gamma B]$. Hence we have $A\cap B=(A\Gamma B]$.\\$\Longleftarrow$. 
Let $A$ be an ideal of $M$. By hypothesis, we have $A=(A\Gamma A]$, 
so $A$ is idempotent. $\hfill\Box$\medskip

\noindent{\bf Theorem 10.} {\it Let M be a $po$-$\Gamma$-semigroup. 
The ideals of M are weakly prime if and only if they are idempotent 
and they form a chain}.\medskip

\noindent{\bf Proof.} $\Longrightarrow$. Let $A$, $B$ be ideals of 
$M$. One can easily prove that $(A\Gamma B]$ is an ideal of $M$.
Since $A,B,(A\Gamma B]$ are ideals of $M$, $A\Gamma B\subseteq 
(A\Gamma B]$ and $(A\Gamma B]$ is weakly prime, we have $A\subseteq 
(A\Gamma B]\subseteq (M\Gamma B]\subseteq (B]=B$ or $B\subseteq 
(A\Gamma B]\subseteq (A\Gamma M]\subseteq (A]=A$, thus the ideals of 
$M$ form a chain. Furthermore, since $A$ and $(A\Gamma A]$ are ideals 
of $M$, $A\Gamma A\subseteq (A\Gamma A]$ and $(A\Gamma A]$ is weakly 
prime, we have $A\subseteq (A\Gamma A]\subseteq (M\Gamma A]\subseteq 
(A]=A,$ so $A=(A\Gamma A]$.\\$\Longleftarrow$. Let $A,B,T$ be ideals 
of $M$ such that $A\Gamma B\subseteq T$. If $A\subseteq B$ then, by 
Lemma 9, $A=A\cap B=(A\Gamma B]\subseteq (T]=T$. If $B\subseteq A$, 
then $B=A\cap B=(A\Gamma B]\subseteq (T]=T$. $\hfill\Box$\medskip

\noindent{\bf Lemma 11.} {\it Let M be a $po$-$\Gamma$-semigroup. If 
M is intra-regular, then$$I(x)=(M\Gamma x\Gamma M] \mbox { for every 
} x\in M.$$}{\bf Proof.} Let $x\in M$. Since $(M\Gamma x\Gamma M]$ is 
an ideal of $M$, by Lemma 5, it is semiprime. Take an element 
$\gamma\in\Gamma$ $(\Gamma\not=\emptyset)$. Since $x\gamma x\in M$ 
and $(x\gamma x)\gamma (x\gamma x)\in (M\Gamma x\Gamma x]$, we have 
$x\gamma x\in (M\Gamma x\Gamma x]$, then $x\in (M\Gamma x\Gamma x]$, 
and $I(x)\subseteq (M\Gamma x\Gamma x]$. On the other hand, $(M\Gamma 
x\Gamma x]\subseteq I(x)$, thus we get $I(x)=(M\Gamma x\Gamma M]$. 
$\hfill\Box$\medskip

\noindent{\bf Lemma 12.} {\it If M is an $po$-$\Gamma$-semigroup, 
$x,y\in M$ and $\gamma\in\Gamma$, then$$I(x\gamma y)\subseteq 
I(x)\cap I(y).$$In particular, if M is intra-regular, then $I(x\gamma 
y)=I(x)\cap I(y)$}.\bigskip

\noindent{\bf Proof.} Since $x\gamma y\in I(x)\Gamma M\subseteq I(x)$ 
and $x\gamma y\in M\Gamma I(y)\subseteq I(y)$, we have $I(x\gamma 
y)\subseteq I(x)\cap I(y)$. Let now $M$ be intra-regular and $t\in 
I(x)\cap I(y)$. Then, by Lemma 11, we have $t\in (M\Gamma x\Gamma M]$ 
and $t\in (M\Gamma y\Gamma M]$. Thus we have$$t\le a\mu x\rho b \mbox 
{ and } t\le c\xi y\zeta d \mbox { for some } a,b,c,d\in M,\, 
\mu,\rho,\xi,\zeta\in\Gamma.$$Then $t\gamma t\le (c\xi y\zeta 
d)\gamma (a\mu x\rho b)=c\xi (y\zeta d\gamma a\mu x)\rho b.$ In 
addition, we have $y\zeta d\gamma a\mu x\in I(x\gamma y)$. Indeed, by 
Lemma 11,$$(y\zeta d\gamma a\mu x)\gamma (y\zeta d\gamma a\mu x)\in 
M\Gamma (x\gamma y)\Gamma M\subseteq {\Big(}M\Gamma (x\gamma y)\Gamma 
M{\Big]}=I(x\gamma y).$$Since $M$ is intra-regular and $I(x\gamma y)$ 
is an ideal of $M$, by Lemma 5, $I(x\gamma y)$ is semiprime. So we 
get $y\zeta d\gamma a\mu x\in I(x\gamma y)$. Since $I(x\gamma y)$ is 
an ideal of $M$, we have $c\xi (y\zeta d\gamma a\mu x)\rho b\in 
M\Gamma I(x\gamma y)\Gamma M\subseteq I(x\gamma y)\Gamma M\subseteq 
I(x\gamma y)$, then $t\gamma t\in I(x\gamma y)$. Since $I(x\gamma y)$ 
is semiprime, we have $t\in I(x\gamma y)$. Thus we get $I(x)\cap 
I(y)\subseteq I(x\gamma y)$ and so $I(x\gamma y)=I(x)\cap I(y)$. 
$\hfill\Box$\medskip

\noindent{\bf Theorem 13.} {\it Let $M$ be a $po$-$\Gamma$-semigroup. 
The ideals of M are prime if and only if  they form a chain and M is 
intra-regular}.  \medskip

\noindent{\bf Proof.} $\Longrightarrow$. The ideals of $M$ are prime, 
so they are weakly prime and semiprime. Since they are weakly prime, 
by Theorem 10, they form a chain. Let now $a\in M$ and 
$\gamma\in\Gamma$. Since $(M\Gamma a\gamma a\Gamma M]$ is an ideal of 
$M$, $a\gamma a\in (M\Gamma a\gamma a\Gamma M]$ and
$(M\Gamma a\gamma a\Gamma M]$ is semiprime, we have $a\gamma a\in 
(M\Gamma a\gamma a\Gamma M]$, and $a\in (M\Gamma a\gamma a\Gamma M]$. 
Thus $M$ is intra-regular. \\
\noindent $\Longleftarrow$. Suppose $M$ is intra-regular and the 
ideals of $M$ form a chain. Let now $T$ be an ideal of $M$, $a,b\in 
M$ and $\gamma\in\Gamma$ such that $a\gamma b\in T$. We have 
$I(a)\subseteq I(b)$ or $I(b)\subseteq I(a)$. If $I(a)\subseteq I(b)$ 
then, by Lemma 12, we have $a\in I(a)=I(a)\cap I(b)=I(a\gamma 
b)\subseteq I(T)=T$. If $I(b)\subseteq I(a)$, then $b\in 
I(b)=I(a)\cap I(b)=I(a\gamma b)\subseteq T$.$\hfill\Box$\medskip

\noindent{\bf Proposition 14.} {\it Let M be an intra-regular 
$po$-$\Gamma$-semigroup. If the ideals of M form a chain then, for 
every $x,y\in M$ and every $\gamma\in\Gamma$, we have$$x\in (M\Gamma 
x\gamma y\Gamma M] \mbox { or } y\in (M\Gamma x\gamma y\Gamma 
M].$$}{\bf Proof.} Assuming the ideals of $M$ form a chain, let
$x,y\in M$ and $\gamma\in\Gamma$. By Theorem 13, the ideals of $M$ 
are prime. Since $(M\Gamma x\gamma y\Gamma M]$ is an ideal of $M$, 
$(M\Gamma x\gamma y\Gamma M]$ is prime. Since $(x\gamma x)\gamma 
(y\gamma y)\in (M\Gamma x\gamma y\Gamma M]$, we have $x\gamma x\in 
(M\Gamma x\gamma y\Gamma M]$ or $y\gamma y\in (M\Gamma x\gamma 
y\Gamma M]$. If $x\gamma x\in (M\Gamma x\gamma y\Gamma M]$ then, 
since $(M\Gamma x\gamma y\Gamma M]$ is prime, we have $x\in (M\Gamma 
x\gamma y\Gamma M]$. If $y\gamma y\in (M\Gamma x\gamma y\Gamma M]$, 
then $y\in (M\Gamma x\gamma y\Gamma M]$.$\hfill\Box$ \medskip

\noindent{\bf Definition 15.} A $po$-$\Gamma$-semigroup $M$ is called 
a {\it chain of simple semigroups} if there exists a semilattice 
congruence $\sigma$ on $M$ such that $(x)_\sigma$ is a simple 
subsemigroup of $M$ for every $x\in M$, and for any $x,y\in M$ and 
any $\gamma\in\Gamma$ we have $$(x)_{\sigma}=(x\gamma y)_{\sigma} 
\mbox { or } (y)_{\sigma}=(x\gamma y)_{\sigma}$$(in other words, the 
set $M/\sigma$ endowed with the relation $$(x)_{\sigma}\preceq 
(y)_{\sigma} \Leftrightarrow (x)_{\sigma}=(x\gamma y)_{\sigma} \; 
\forall \; \gamma\in\Gamma$$ is a chain).\medskip

\noindent{\bf Theorem 16.} {\it A $po$-$\Gamma$-semigroup M is 
intra-regular and the ideals of M form a chain if and only if M is 
chain of simple semigroups}.\medskip

\noindent{\bf Proof.} $\Longrightarrow$. Since $M$ is intra-regular 
and $\cal N$ is a semilattice congruence on $M$, by Theorem 
$8(1)\Rightarrow (5)$, $(x)_{\cal N}$ is a simple subsemigroup of $M$ 
for every $x\in M$, so $M$ is a semilattice of simple semigroups. Let 
now $x,y\in M$ and $\gamma\in\Gamma$. By Proposition 14, we have 
$x\in (M\Gamma x\gamma y\Gamma M]$ or $y\in (M\Gamma x\gamma y\Gamma 
M]$. If $x\in (M\Gamma x\gamma y\Gamma M]$, then $$N(x)\ni x\le a\mu 
(x\gamma y)\rho b \mbox { for some } a,b\in M, \mu,\,\rho\in\Gamma.$$ 
Since $N(x)$ is a filter of $M$, we have $a\mu(x\gamma y)\rho b\in 
N(x)$, $x\gamma y\in N(x)$ and $N(x\gamma y)\subseteq N(x)$. If $y\in 
(M\Gamma x\gamma y\Gamma M]$, then $$N(y)\ni y\le c\xi (x\gamma 
y)\zeta d \mbox { for some } c,d\in M, \xi,\zeta\in\Gamma,$$ then 
$c\xi (x\gamma y)\zeta d\in N(y)$, $x\gamma y\in N(y)$, and 
$N(x\gamma y)\subseteq N(y)$. On the other hand, since $x\gamma y\in 
N(x\gamma y)$, we have $x\in N(x\gamma y)$ and $y\in N(x\gamma y)$, 
so $N(x)\subseteq N(x\gamma y)$ and $N(y)\subseteq N(x\gamma y)$. 
Thus we get $N(x)=N(x\gamma y)$ or  $N(y)=N(x\gamma y)$, then 
$(x)_{\cal N}=(x\gamma y)_{\cal N}$ or $(y)_{\cal N}=(x\gamma 
y)_{\cal N}$. Therefore $M$ is a chain of simple 
semigroups.\smallskip

\noindent $\Longleftarrow$. Let $\sigma$ be a semilattice congruence 
on $M$ such that $(x)_{\sigma}$ is a simple subsemigroup of $M$ for 
every $x\in M$ and $(M/\sigma,\preceq)$ is a chain. By Theorem 13 it 
is enough to prove that the ideals of $M$ are prime. So, let $I$ be 
an ideal of $M$, $a,b\in M$ and $\gamma\in\Gamma$ such that $a\gamma 
b\in I$. The set $(a\gamma b)_{\sigma}\cap I$ is an ideal of 
$(a\gamma b)_{\sigma}$. Indeed:

\hspace{0.3cm}$\emptyset\not=(a\gamma b)_{\sigma}\cap I\subseteq 
(a\gamma b)_{\sigma}$ ($a\gamma b\in (a\gamma b)_{\sigma}$, $a\gamma 
b\in I)$
\begin{eqnarray*}{\Big(}(a\gamma b)_{\sigma}\cap I{\Big)}\Gamma 
(a\gamma b)_{\sigma}&\subseteq&(a\gamma b)_{\sigma}\Gamma (a\gamma 
b)_{\sigma}\cap I\Gamma (a\gamma b)_{\sigma}=(a\gamma b)_{\sigma}\cap 
I\Gamma (a\gamma b)_{\sigma}\\&\subseteq& (a\gamma b)_{\sigma}\cap 
I\Gamma M\subseteq (a\gamma b)_{\sigma}\cap I,
\end{eqnarray*}\begin{eqnarray*}(a\gamma b)_{\sigma}\Gamma 
{\Big(}(a\gamma b)_{\sigma}\cap I{\Big)}&\subseteq& (a\gamma 
b)_{\sigma}\Gamma (a\gamma b)_{\sigma}\cap (a\gamma b)_{\sigma}\Gamma 
I\subseteq (a\gamma b)_{\sigma}\cap M\Gamma I\\&\subseteq& (a\gamma 
b)_{\sigma}\cap I.\end{eqnarray*}Let $x\in (a\gamma b)_{\sigma}\cap 
I$ and $(a\gamma b)_{\sigma}\ni y\le x$. Since $M\ni y\le x\in I$ and 
$I$ is an ideal of $M$, we have $y\in I$, then $y\in (a\gamma 
b)_{\sigma}\cap I$. Since $(a\gamma b)_{\sigma}$ is simple, we have 
$(a\gamma b)_{\sigma}\cap I=(a\gamma b)_{\sigma}$. By hypothesis, 
$(a)_{\sigma}=(a\gamma b)_{\sigma}$ or $(b)_{\sigma}=(a\gamma 
b)_{\sigma}$. Then we have $a\in I$ or $b\in I$, thus $I$ is a prime. 
$\hfill\Box$ \medskip

\noindent{\bf Lemma 17.} {\it Let M be an $po$-$\Gamma$-semigroup, T 
a subsemigroup of M and $x\in T$. Then the set $(M\Gamma x\Gamma 
M]\cap T$ is an ideal of T}.\medskip

\noindent{\bf Proof.} The set $(M\Gamma x\Gamma M]\cap T$ is a 
nonempty subset of $M$. This is because, for any $\gamma\in\Gamma$, 
$x\gamma x\gamma x\in M\Gamma x\Gamma M$ and $x\gamma x\gamma x\in 
(T\Gamma T)\Gamma T\subseteq T\Gamma T\subseteq T$. Moreover,
\begin{eqnarray*}{\Big(}(M\Gamma x\Gamma M]\cap T{\Big)}\Gamma 
T&\subseteq&(M\Gamma x\Gamma M]\Gamma T\cap T\Gamma T\subseteq 
(M\Gamma x\Gamma M]\Gamma (M]\cap T\\&\subseteq& {\Big(}M\Gamma 
x\Gamma (M\Gamma M){\Big]}\cap T\subseteq (M\Gamma x\Gamma M]\cap T,
\end{eqnarray*}In a similar way, we have $T\Gamma {\Big(}(M\Gamma 
x\Gamma M]\cap T{\Big)}\subseteq (M\Gamma x\Gamma M]\cap T.$ Let now 
$a\in (M\Gamma x\Gamma M]\cap T$ and $T\ni b\le a$. Since $a\in 
(M\Gamma x\Gamma M]$, there exist $u,v\in M$ and $\xi,\zeta\in\Gamma$ 
such that $a\le u\xi x\zeta v$. Then we have $b\le u\xi x\zeta v\in 
M\Gamma x\Gamma M$, and $b\in (M\Gamma x\Gamma 
M]$.$\hfill\Box$\medskip

\noindent{\bf Theorem 18.} {\it Let M be an intra-regular 
$po$-$\Gamma$-semigroup. Then the set $(x)_{\cal N}$ is a maximal 
simple subsemigroup of M for every $x\in M$. Conversely, if T is a 
maximal simple subsemigroup of M, then there exists $x\in M$ such 
that $T=(x)_{\cal N}$.}\medskip

\noindent{\bf Proof.} $\Longrightarrow$. Let $x\in M$. By the Theorem 
$8(1)\Rightarrow (5)$, the set $(x)_{\cal N}$ is a simple 
subsemigroup of $M$. Let now $T$ be a simple subsemigroup of $M$ such 
that $T\supseteq (x)_{\cal N}$. Then $T=(x)_{\cal N}$. Indeed: Let 
$y\in T$. Since $x\in T$, by Lemma 17, the set $(M\Gamma x\Gamma 
M]\cap T$ is an ideal of $T$. Since $T$ is a simple, we have 
$(M\Gamma x\Gamma M]\cap T=T$, then $y\in (M\Gamma x\Gamma M]$. Since 
$M$ is intra-regular, by Lemma 3, we have $x\in N(y)$, and 
$N(x)\subseteq N(y)$. On the other hand, since $y\in T$, the set 
$(M\Gamma y\Gamma M]\cap T$ is an ideal of $T$. So $(M\Gamma y\Gamma 
M]\cap T=T$, $x\in (M\Gamma y\Gamma M]$, $y\in N(x)$, and 
$N(y)\subseteq N(x)$. Therefore we have $N(x)=N(y)$, so $y\in 
(x)_{\cal N}$, and $T\subseteq (x)_{\cal N}$. Since $(x)_{\cal N}$ is 
a subsemigroup of $M$, we get $T=(x)_{\cal N}$.\smallskip

\noindent$\Longleftarrow$. Let $T$ be a maximal simple subsemigroup 
of $M$. Take an element $x\in T$ $(T\not=\emptyset)$. Exactly as in 
the proof of the ``$\Rightarrow$"-part of the theorem given above, we 
prove that $T\subseteq (x)_{\cal N}$. Since $(x)_{\cal N}$ is a 
simple subsemigroup of $M$ (cf. Theorem $8(1)\Rightarrow (5)$) and 
$T$ is a maximal simple subsemigroup of $M$, we have $T=(x)_{\cal 
N}$. $\hfill\Box$\\
\noindent{\bf Corollary 19.} {\it For an intra-regular 
$po$-$\Gamma$-semigroup M, the set $\{(x)_{\cal N} \mid x\in M\}$ 
coincides with the set of all maximal simple subsemigroup of 
M.}\medskip

\noindent{\bf Definition 20.} [6] A $\Gamma$-semigroup $M$ is called 
{\it left} (resp. {\it right}) {\it regular} if $$x\in (M\Gamma 
x\gamma x] \mbox { (resp. } x\in (x\gamma x\Gamma M])$$for every  
$x\in M$ and every $\gamma\in\Gamma$.\medskip

\noindent In a similar way as in Theorem 8 above and the Theorem 6 in 
[6], we can prove the following \medskip

\noindent{\bf Theorem 21.} {\it Let M be a $po$-$\Gamma$-semigroup. 
The following are equivalent:

$(1)$ M is left regular and left duo.

$(2)$ $N(x)=\{y\in M \mid x\in (M\Gamma y]\}$ for every
$x\in M$.

$(3)$ $\cal N=\cal L$.

$(4)$ For every left ideal L of M, we have $L=
\bigcup\limits_{x \in L} {(x)_{\cal N}}$.

$(5)$ $(x)_{\cal N}$ is a left simple subsemigroup of M for
every $x\in M$.

$(6)$ $M$ is a semilattice of left simple semigroups.

$(7)$ Every left ideal of M is semiprime and two-sided.}\medskip

\noindent The right analogue of the above theorem also holds.
{\small
\bigskip

\medskip

\noindent University of Athens, Department of Mathematics, 15784 
Panepistimiopolis, Greece\\
email: nkehayop@math.uoa.gr

\end{document}